\documentclass[12pt,oneside,english]{amsart}
\usepackage[latin1]{inputenc}
\usepackage{amssymb}

\makeatletter
\newtheorem{theorem}{Theorem}
\newtheorem{lemma}{Lemma}
\newtheorem{algorithm}{Algorithm}
\newtheorem{example}{Example}

\numberwithin{equation}{section}
\usepackage{babel}

\makeatother
\begin{document}
\baselineskip=17pt

\title[Generalizations of the Rowland theorem]{Generalizations of the Rowland theorem}

\author{Vladimir Shevelev}
\address{Departments of Mathematics \\Ben-Gurion University of the
 Negev\\Beer-Sheva 84105, Israel. e-mail:shevelev@bgu.ac.il}


\begin{abstract}
We prove the theorems which are equivalent to the Roland's results such that a new form of them allows to consider some generalizations. In particular, we give generators of primes more than a fixed prime.
\end{abstract}

\maketitle

\section{Introduction }

In 2008, Rowland [1], using very elementary tools, discovered a very interesting fact.
\begin{theorem}\label{1}
Let $a(1)=7$ and, for $n\geq2,$
$$a(n)=a(n-1)+\gcd (n, \enskip a(n-1)). $$
Then, for every $n\geq2,$ the difference $a(n)-a(n-1)$ is 1 or prime.
 \end{theorem}
 He also mentioned that, in particular, the following similar theorem is true (it also follows from his proof).
 \begin{theorem}\label{2}
Let $b(3)=6$ and, for $n\geq4,$
$$b(n)=b(n-1)+\gcd (n, \enskip b(n-1)). $$
Then, for every $n\geq2,$ the difference $b(n)-b(n-1)$ is 1 or prime.
 \end{theorem}
 In spite of the nearness of Theorem 1 and 2, there is an essential distinction: $\lim\sup a(n)/n=3,$ while $\lim\sup b(n)/n=2.$
 An infinite sets of the initial conditions for which one of Theorems 1,2 is true we considered in [2].\newline
 \indent \slshape A generalization of Theorems 1,2 could $appear$ if $to$ become free from the hard dependence on the value of \upshape $\gcd.$\newline
 \indent It is important that the Rowland's proof allows to write his results in such a form. The following theorems are equivalent to Theorems 1,2.\newline
  \bfseries Theorem 1a. \mdseries \slshape Let $a(1)=7$ and, for $n\geq2,$
$$a(n)=\begin{cases}a(n-1)+1,\enskip if \;\;\gcd (n, a(n-1))=1
\\3n,\enskip \;\; otherwise\end{cases}.  $$
Then, for every $n\geq2,$ the difference $a(n)-a(n-1)$ is 1 or prime.\upshape\newline
\bfseries Theorem 2a. \mdseries \slshape Let $b(3)=6$ and, for $n\geq4,$
$$b(n)=\begin{cases}b(n-1)+1,\enskip if \;\;\gcd (n, b(n-1))=1
\\2n,\enskip \;\; otherwise\end{cases}.  $$
Then, for every $n\geq4,$ the difference $b(n)-b(n-1)$ is 1 or prime.\upshape
\section{Direct proof of Theorems $1a, \enskip2a$ with general initials}
Excepting a trivial case when $a(n)=1$ identically, the sequences which are defined by Theorems 1a and 2a contain the values in which $a(n)=3n$ or $a(n)=2n$ correspondingly. Therefore the general initial conditions for them have the form $a(n_1)=3n_1$ or $a(n_1)=2n_1$ correspondingly. Below we give a direct and simple proofs of Theorems $1a, \enskip2a$ with general initial conditions. Our proof is again based on Rowland idea, but in some another form.\newline
 \bfseries Theorem 1b. \mdseries \slshape Let, for $n_1\geq2, \enskip a(n_1)=3n_1$ and, for $n\geq n_1+1,$
 \begin{equation}\label{2.1}
a(n)=\begin{cases}a(n-1)+1,\enskip if \;\;\gcd (n, a(n-1))=1,
\\3n,\enskip \;\; otherwise\end{cases}.
\end{equation}
Then, for every $n\geq n_1+1,$ the difference $a(n)-a(n-1)$ is 1 or prime.\upshape\newline
\bfseries Theorem 2b. \mdseries \slshape Let, for $n_1\geq3, \enskip b(n_1)=2n_1$ and, for $n\geq n_1+1,$
 \begin{equation}\label{2.2}
b(n)=\begin{cases}b(n-1)+1,\enskip if \;\;\gcd (n, b(n-1))=1,
\\2n,\enskip \;\; otherwise\end{cases}.
\end{equation}
Then, for every $n\geq n_1+1,$ the difference $b(n)-b(n-1)$ is 1 or prime.\upshape\newline
\bfseries Proof of Theorem 2b. \mdseries
 Let $k$ be the smallest positive integer such that
  \begin{equation}\label{2.3}
\gcd(n_1+k, \enskip 2n_1+k-1)=d>1.
\end{equation}
Then from (2.2) and {2.3) we have
$$b(n_1)=2n_1,$$
$$b(n_1+1)=2n_1+1,$$
$$...$$
$$b(n_1+k-1)=2n_1+k-1,$$
\begin{equation}\label{2.4}
b(n_1+k)=2(n_1+k).
\end{equation}
Besides, by (2.3), $\gcd(n_1-1,\enskip n_1+k)=d$ and, therefore,
\begin{equation}\label{2.5}
\gcd(n_1-1,\enskip k+1)=d.
\end{equation}
Thus some prime divisor $P$ of $n_1-1$ divides $k+1,$ i.e. $k+1\geq P.$ All the more,
\begin{equation}\label{2.6}
k+1\geq p,
\end{equation}
where $p$ is the smallest prime divisor of $n_1-1$ (by the condition, $n_1-1\geq2$).
\indent On the other hand, inside (2.4) there is a row
\begin{equation}\label{2.7}
b(n_1+p-2)=2n_1+p-2.
\end{equation}\newpage
Since $\gcd(2n_1+p-2,\enskip n_1+p-1)\geq p,$ then, by the definition of $k,$ we have $k\leq p-1$ and, in view of (2.6),
\begin{equation}\label{2.8}
k+1=p.
\end{equation}
Now from (2.4) we find
$$ b(n_1+k)-b(n_1+k-1)=k+1=p,$$
$$b(n_1+k)=2(n_1+k)$$
and, by the evident induction, we are done. $\blacksquare$\newline
\bfseries Proof of Theorem 1b \mdseries is the same if to replace (2.4) by
 $$a(n_1+k-1)=3n_1+k-1,$$

$$a(n_1+k)=3(n_1+k),$$

(2.5) by

$$\gcd(2n_1-1,\enskip k+1)=d$$

and (2.7) by

$$a(n_1+\frac {p-1}{2}-1)=3n_1+\frac {p-1} {2}-1,$$

 where now $p$ is the smallest prime divisor of $2n_1-1\geq3,$ such that $2k+1=p. \blacksquare$\newline
\section{Algorithms of fast calculation of nontrivial increments in sequences of Theorems $1b,\enskip2b$}
Sequences of the considered type contain too many points of trivial $1$-increments. Nevertheless, it is interesting to see quicker what primes these sequences contain. Therefore, the following problem is actual from the computation point of view : to accelerate this algorithm for receiving of primes by the omitting of the trivial increments. Using induction, from the proofs of theorems $1b,\enskip2b$ we discover a sense of these primes by the following algorithms.\newline
\indent For integer $n\geq2,$ denote $n^*$ the least prime divisor of $n.$ We start with an algorithm of fast calculation of the Rowland's primes which are obtained by Theorem 1.

\begin{algorithm}\label{1}
Put
$$N_1=4,\enskip p_1=5$$\newpage
 and, for $i\geq1,$ put
 $$N_{i+1}=N_i+p_i-1,\enskip p_{i+1}=(N_{i+1}+1)^*.$$
 Then sequence $\{p_i\}_{i\geq1}$ is the sequence of nontrivial increments of the Rowland sequence, which is defined by Theorem 1.
 \end{algorithm}
By this algorithm, we consecutively find
$$p_1=5,\enskip N_2=8, \enskip p_2=9^*=3,\enskip N_3=10,\enskip p_3=11,\enskip N_4=20,\enskip p_4=3,...$$
Using formulas of Algorithm 1, it is easy to obtain the following recursion for $p_n$ with the automatically produced initial term:
\begin{equation}\label{3.1}
 p_n=(6-n+\sum_{i=1}^{n-1}p_i)^*.
\end{equation}
\begin{example}\label{1}
We have consecutively (cf. sequence $A137613$ in $[3]$):
$$p_1 =(6-1)^*=5,$$
$$p_2=(6-2+5)^*=3,$$
$$p_3=(6-3+5+3)^*=11,$$
$$p_4=(6-4+5+3+11)^*=3,$$
$$p_5=(6-5+5+3+11+3)^*=23,$$
etc.
\end{example}
Note that, beginning with $n\geq3,$ the Rowland's sequence could be defined by the initial condition $a(3)=9.$ Therefore, the following algorithm is a generalization of Algorithm 1.
\begin{algorithm}\label{2} Consider sequence of Theorem $1b$ with initial condition $a(n_1)=3n_1.$
Put
$$N_1=2(n_1-1),\enskip p_1=(N_1+1)^*$$
 and, for $i\geq1,$ put
 $$N_{i+1}=N_i+p_i-1,\enskip p_{i+1}=(N_{i+1}+1)^*.$$
 Then sequence $\{p_i\}_{i\geq1}$ is the sequence of nontrivial increments of the sequence, which is defined by Theorem $1b.$
 \end{algorithm}
E.g., for $n_1=4,$ by this algorithm, we consecutively find
$$N_1=6,\enskip p_1=7,\enskip N_2=12, \enskip p_2=13,\enskip N_3=24,\enskip p_3=5,\enskip N_4=28,\enskip p_4=29,...$$
Using formulas of Algorithm 2, we obtain the following recursion for $p_n$ with the automatically produced initial term:\newpage
\begin{equation}\label{3.2}
 p_n=(2n_1-n+\sum_{i=1}^{n-1}p_i)^*,
\end{equation}
such that, for $n_1=3, \enskip n\geq3,$ we have (3.1).
\begin{example}\label{2}
In case of $n_1=4,$ we have consecutively:
$$p_1 =(8-1)^*=7,$$
$$p_2=(8-2+7)^*=13,$$
$$p_3=(8-3+7+13)^*=5,$$
$$p_4=(8-4+7+13+5)^*=29,$$
etc.
\end{example}
Finally, we give the corresponding algorithm for sequences which are defined by Theorem $2b.$
\begin{algorithm}\label{3} Consider sequence of Theorem $2b$ with initial condition $a(n_1)=2n_1.$
Put
$$N_1=n_1-2,\enskip p_1=(N_1+1)^*$$
 and, for $i\geq1,$ put
 $$N_{i+1}=N_i+p_i-1,\enskip p_{i+1}=(N_{i+1}+1)^*.$$
 Then sequence $\{p_i\}_{i\geq1}$ is the sequence of nontrivial increments of the sequence, which is defined by Theorem $2b.$
 \end{algorithm}
E.g., for $n_1=5,$ by this algorithm, we consecutively find
$$N_1=3,\enskip p_1=2,\enskip N_2=4, \enskip p_2=5,\enskip N_3=8,\enskip p_3=3,\enskip N_4=10,\enskip p_4=11,...$$
Here we obtain the following recursion for $p_n:$
\begin{equation}\label{3.3}
 p_n=(n_1-n+\sum_{i=1}^{n-1}p_i)^*.
\end{equation}
\begin{example}\label{3}
In case of $n_1=5,$ we have consecutively:
$$p_1 =(5-1)^*=2,$$
$$p_2=(5-2+2)^*=5,$$
$$p_3=(5-3+2+5)^*=3,$$
$$p_4=(5-4+2+5+3)^*=11,$$
etc.
\end{example}
\newpage
 \section{$\nu$- Generalizations: other generators of primes }
 \begin{theorem}\label{3} Let $\nu$ be a positive even integer. Let, for $n_1\geq2, \enskip c(n_1)=3n_1+\nu$ and, for $n\geq n_1+1,$
 \begin{equation}\label{4.1}
c(n)=\begin{cases}3n+\nu,\enskip if \;\;\gcd (n,\enskip c(n-1))>1 \enskip and\enskip \rho\geq\frac{\nu-2} {2},
\\c(n-1)+1,\enskip \;\; otherwise\end{cases}
\end{equation}
where, for $n>n_1,$
\begin{equation}\label{4.2}
 \rho=\rho(n):=n-\max\{l<n: \enskip c(l)=3l+\nu\}.
 \end{equation}
 Then, for every $n\geq n_1+1,$ the difference $c(n)-c(n-1)-\nu$ is $1-\nu$ or prime.
 \end{theorem}
 Note that, condition (4.2) means that in this sequence the distance between the consecutive nontrivial increments (i.e. the increments which are different from 1) is not less than $\frac{\nu-2} {2}.$\newline
 \bfseries Proof. \mdseries Let $k$ be the smallest positive integer not less then $\frac{\nu-2} {2}$ such that
  \begin{equation}\label{4.3}
\gcd(n_1+k, \enskip 3n_1+k-1+\nu)=d>1.
\end{equation}
Then from (4.1) and {4.3) we have
$$c(n_1)=3n_1+\nu,$$
$$c(n_1+1)=3n_1+1+\nu,$$
$$...$$
$$c(n_1+k-1)=3n_1+k-1+\nu,$$
\begin{equation}\label{4.4}
c(n_1+k)=3(n_1+k)+\nu.
\end{equation}
 Besides, by (4.3), $\gcd(2n_1+\nu-1,\enskip n_1+k)=d>1$ and, therefore,

\begin{equation}\label{4.5}
\gcd(2n_1+\nu-1,\enskip 2k-\nu+1)=d>1.
\end{equation}
Thus some prime divisor $P$ of $2n_1+\nu-1$ divides $2k-\nu+1,$ i.e. $2k-\nu+1\geq P.$ All the more,
\begin{equation}\label{4.6}
2k-\nu+1\geq p,
\end{equation}
where $p$ is the smallest prime divisor of $2n_1+\nu-1$ (by the condition, $\nu$ is even, therefore $p$ is odd).
Note that (4.6) shows that \slshape the condition $k\geq\frac {\nu-2} {2}$ is necessary.\upshape\newline
\indent On the other hand, inside (2.4) there is a row
\begin{equation}\label{4.7}
c(n_1+\frac {p+\nu-3}{2})=3n_1+\frac {p+\nu-3}{2}+\nu=\frac{1} {2}(6n_1+3\nu-3+p).
\end{equation}
 Since $n_1+\frac {p+\nu-1}{2}=\frac {1}{2}(2n_1+\nu-1+p),$ then we see that \newpage
 $$\gcd(c(n_1+\frac {p+\nu-3}{2}), \enskip n_1+\frac{p+\nu-1}{2})\geq p,$$
 and, by the definition of $k,$ we have $k\leq\frac{p+\nu-1} {2},$ and, in view of (2.6), we conclude that
\begin{equation}\label{4.8}
2k-\nu+1=p.
\end{equation}
Now from (2.4) we find
$$ c(n_1+k)-c(n_1+k-1)-\nu=2k+1-\nu=p,$$
and, for $n_2:=n_1+k,$
$$c(n_2)=3(n_2)+\nu$$
and, by the evident induction, we are done. $\blacksquare$\newline
 \indent The following theorem is proved quite analogously.
  \begin{theorem}\label{4} Let $\nu$ be a positive integer. Let, for $n_1\geq2, \enskip c(n_1)=2n_1+\nu$ and, for $n\geq n_1+1,$
 \begin{equation}\label{4.9}
c(n)=\begin{cases}2n+\nu,\enskip if \;\;\gcd (n,\enskip c(n-1))>1 \enskip and\enskip \rho\geq\nu+1,
\\c(n-1)+1,\enskip \;\; otherwise\end{cases}
\end{equation}
where, for $n>n_1,$
\begin{equation}\label{4.10}
 \rho=\rho(n):=n-\max\{l<n: \enskip c(l)=2l+\nu\}.
 \end{equation}
 Then, for every $n\geq n_1+1,$ the difference $c(n)-c(n-1)-\nu$ is $1-\nu$ or prime.
 \end{theorem}
 \section{$\nu$- Generalizations: generators of primes more than $p_m$}
 The following theorem gives generators of primes $p>p_m,\enskip m\geq3,$ where $p_n$ is the $n$-th prime.
  \begin{theorem}\label{5} Let $m\geq3,\enskip\nu\geq0$ be even and not exceeding $ p_m-3.$  Let, for $n_1>\nu+2, \enskip c(n_1)=3n_1-\nu$ and, for $n\geq n_1+1,$
 \begin{equation}\label{5.1}
c(n)=\begin{cases}3n-\nu,\enskip if \;\;\gcd (n,\enskip c(n-1))>1 \enskip and\enskip \gcd(n,\enskip\prod_{i=1}^{m}p_i )=1,
\\c(n-1)+1.\enskip \;\; otherwise\end{cases}
\end{equation}
Then, for every $n\geq n_1+1,$ the enlarged on $\nu$ difference $c(n)-c(n-1)$ is $1+\nu$ or prime more than $p_m.$
 \end{theorem}
 \bfseries Proof \mdseries is similar to proof of Theorem 3 with the replacing $\nu$ by $-\nu.$
  Let $k$ be the smallest positive integer such that

$$\gcd(n_1+k, \enskip \prod_{i=1}^{m}p_i )=1,  \enskip and \enskip$$\newpage
\begin{equation}\label{5.2}
\gcd(n_1+k, \enskip 3n_1+k-1+\nu)=d>1.
\end{equation}
Then, as in proof of Theorem 3, we conclude that

$$\gcd(2n_1-\nu-1,\enskip 2k+\nu+1)=d>1  \enskip and \enskip$$
\begin{equation}\label{5.3}
\gcd(2n_1-\nu-1, \enskip \prod_{i=1}^{m}p_i )=1.  \enskip
\end{equation}
Thus if $p$ is the smallest prime divisor of $2n_1-\nu-1,$ then
\begin{equation}\label{5.4}
p>p_m\geq \nu+3,
\end{equation}
and
\begin{equation}\label{5.5}
2k+\nu+1\geq p.
\end{equation}
Note that the condition $n_1>\nu+2$ is necessary. Indeed, by (5.3), we have  $2n_1-\nu-1\geq p>p_m\geq\nu+3$ and the inequality $n_1>\nu+2$ follows.\newline
\indent On the other hand, \slshape in view of (5.4)\upshape \enskip inside (2.4) with the replacing $\nu$ by $-\nu$ we find a row
\begin{equation}\label{5.6}
c(n_1+\frac {p-\nu-3}{2})=3n_1+\frac {p-\nu-3}{2}-\nu=\frac{1} {2}(6n_1-3\nu-3+p).
\end{equation}
 Since $n_1+\frac {p-\nu-1}{2}=\frac {1}{2}(2n_1-\nu-1+p),$ then we see that $$\gcd(c(n_1+\frac {p-\nu-3}{2}), \enskip n_1+\frac{p-\nu-1}{2})\geq p,$$ and, by the definition of $k,$ we have $k\leq\frac{p-\nu-1} {2},$ and, in view of (5.5), we conclude that
\begin{equation}\label{5.7}
2k+\nu+1=p.
\end{equation}
Now we again obtain that
$$ c(n_1+k)-c(n_1+k-1)+\nu=2k+1+\nu=p,$$
and, for $n_2:=n_1+k,$
$$c(n_2)=3(n_2)-\nu$$
and, by the evident induction, we are done. $\blacksquare$\newline
\indent The following theorem is proved quite analogously.
\begin{theorem}\label{6}
 Let $m\geq2,\enskip\nu$ be even not exceeding $ p_m-2.$  Let, for $n_1\geq2\nu+4, \enskip c(n_1)=2n_1-\nu$ and, for $n\geq n_1+1,$
 \begin{equation}\label{5.8}
c(n)=\begin{cases}2n-\nu,\enskip if \;\;\gcd (n,\enskip c(n-1))>1 \enskip and\enskip \gcd(n,\enskip\prod_{i=1}^{m}p_i )=1,
\\c(n-1)+1.\enskip \;\; otherwise\end{cases}
\end{equation}\newpage
Then, for every $n\geq n_1+1,$ the enlarged on $\nu$ difference $c(n)-c(n-1)$ is $1+\nu$ or prime more than $p_m.$
 \end{theorem}
 Note that the conditions $\nu\leq p_m-2$ and $n_1\geq2\nu+4$ are necessary (cf. proof of Theorem 5).

\section{Examples}
1)An illustration to Theorem 3. $\nu=4,\enskip n_1=6.$ The first terms of the sequence are:
$$ 22,23,24,31,32,33,34,35,46,...$$
with the first nontrivial diminished on 4 increments $3,7,31,...$\newline
2)An illustration to Theorem 4. $\nu=5,\enskip n_1=15.$ The first terms of the sequence are:
$$35,36,37,38,39,...$$
Note that $a(19)=39,$ since $19-15<\nu+1=6.$\newline
The first nontrivial diminished on 5 increments are $19,3,7,3,67,...$\newline
3)An illustration to Theorem 5 (cf. sequences A168143-A168144 in [3]). $p_7=17, \enskip \nu=14,\enskip n_1=17.$ The first terms of the sequence are:
$$37,38,43,44,45,46,55,56,57,58,73,...$$
with the first nontrivial enlarged on $14$ increments $19,23,29,43,71,...$\newline
\indent \slshape Astonishingly that, probably, it is difficult to find a regularity in order to obtain a $\nu $-generalization of algorithms of Section 3.\upshape

\section{An additional result}
An attempt to construct of similar sequences ( in style of Theorems 3-6)  with \slshape increasing\upshape \enskip nontrivial increments leads us to a surprising thing. More exactly, consider the following sequence. Let $a(3)=6$ and, for $n\geq4,$
\begin{equation}\label{7.1}
c(n)=2n,\enskip if \enskip\gcd (n,\enskip c(n-1))>1
\end{equation}
 and for $ m<n,$
  $$ \gcd (m,\enskip c(m-1))<\gcd (n,\enskip c(n-1)),$$ and, otherwise,
$$c(n)=c(n-1)+1.$$
\begin{theorem}\label{7}
 The nontrivial increments of $\{c(n)\}$ form the sequence
 $$\{2^n+1\}_{n\geq1}.$$
\end{theorem}
Firstly, we prove a lemma.\newpage
\begin{lemma}
For $j=1,2,...,n-1,$ we have
$$ \gcd(n+j,\enskip 2n+j-1)\leq n/2.$$
\end{lemma}
\bfseries Proof. \mdseries Indeed, evidently
$$\gcd(n+j,\enskip 2n+j-1)=\gcd(n+j,\enskip n-1)=\gcd(j+1,\enskip n-1),$$
and the lemma follows.$\blacksquare$\newline
\bfseries Proof of Theorem 7.\mdseries  \enskip We use induction with the base $c(3)=6,\enskip c(4)=8.$  We see that to obtain $c(4)=8$ from $c(3)=6$ we add 1 to the value of argument 3 and 0 to to the term 6 of the sequence; after that we obtain $c(4)=c(3)+\gcd (3+1,\enskip c(3)+0)=6+\gcd(4,6)=8.$ Suppose that in case of $c(n)=2n$ we add $n-2$ to the value of argument $n$ and $n-3$ to the term $c(n).$ Then we have $$c(2n-2)=3n-3+\gcd(2n-2,\enskip3n-3)=4n-4.$$
By the same way, putting $2n-2=m,$ we have
 $$c(2m-2)=3m-3+\gcd(2m-2,\enskip3m-3)=4m-4.$$
 Using Lemma 1, we see that
 $$ \gcd(m+j,\enskip 2m+j-1)<m/2=n-1,$$
 i.e. the value of $\gcd$ of the precede step.
 Now, by the induction, we conclude that the numbers $2,\enskip 2\cdot2-1=3,\enskip 2\cdot3-1=5,...,2^k+1,... $ are the smallest values of $\gcd(n,\enskip c(n-1),$ to which do not precede any non-smaller values of $\gcd.\blacksquare$

Close results we could prove, putting, e.g., $c(7)=12$ (here we obtain numbers $4,7,13,25,...$) or $c(5)=12$ with the replacing in (6.1) $c(n)=2n$ by $c(n)=3n$ (here we obtain numbers $6,11,21,41,...$) etc.
\;\;\;\;\;\;\;\;

\;\;\;\;\;\;\;\;


\begin{thebibliography}{1}
\bibitem 1. E. S. Rowland, \slshape\enskip A natural prime-generating recurrence \upshape J.Integer Seq., v.11(2008), Article 08.2.8
\bibitem 2.  V. Shevelev, \slshape An infinite set of generators of primes based on the Roland idea and
conjectures concerning twin primes,\upshape\enskip http://www.arxiv.org/abs/0910.4676 [math. NT].

\bibitem {3}.  N.\enskip J.\enskip A.\enskip Sloane,\enskip\slshape The On-Line Encyclopedia of Integer Sequences \upshape $(http: //www.research.att.com/\sim njas)$
\end{thebibliography}
\end{document}